\documentclass[12pt,reqno]{tran-l}
\usepackage{amssymb}
\usepackage[all]{xy}
\usepackage{times}

\newtheorem*{thm}{Theorem}
\newtheorem*{prop}{Proposition}
\newtheorem*{lem}{Lemma}
\newtheorem*{cor}{Corollary}

\theoremstyle{definition}

\newtheorem*{ack}{Acknowledgments}

\theoremstyle{remark}

\numberwithin{equation}{section}

\newcommand{\Ker}{\operatorname{Ker}}

\newcommand{\Hom}{\operatorname{Hom}}

\newcommand{\ann}{\operatorname{ann}}
\newcommand{\Spec}{\operatorname{Spec}}

\newcommand{\GL}{\operatorname{GL}}

\newcommand{\SL}{\operatorname{SL}}

\newcommand{\rank}{\operatorname{rank}}
\newcommand{\height}{\operatorname{height}}
\newcommand{\depth}{\operatorname{depth}}
\newcommand{\diag}{\operatorname{diag}}

\newcommand{\res}{\operatorname{res}}

\newcommand{\ch}{\operatorname{char}}

\newcommand{\into}{\hookrightarrow}

\newcommand{\gen}[1]{\langle{#1}\rangle}

\renewcommand{\bar}[1]{\overline{#1}}

\renewcommand{\phi}{\varphi}
\newcommand{\osum}{\sigma}

\newcommand{\G}{\mathcal{G}}
\newcommand{\cO}{\mathcal{O}}

\newcommand{\cP}{\mathcal{P}}

\renewcommand{\H}{\mathcal{H}}

\newcommand{\cX}{\mathfrak{X}}

\newcommand{\Sy}{\mathcal{S}}
\newcommand{\Sym}{\mathsf{S}}

\newcommand{\R}{\mathcal{R}}
\newcommand{\A}{\mathcal{A}}

\newcommand{\fp}{\mathfrak{p}}
\newcommand{\fP}{\mathfrak{P}}
\newcommand{\fa}{\mathfrak{a}}

\newcommand{\x}{\mathbf{x}}
\renewcommand{\k}{\Bbbk}
\newcommand{\FF}{\mathbb{F}}

\newcommand{\CC}{\mathbb{C}}
\newcommand{\ZZ}{\mathbb{Z}}
\newcommand{\QQ}{\mathbb{Q}}
\newcommand{\RR}{\mathbb{R}}
\newcommand{\HH}{\mathbb{H}}

\newcommand{\kL}{\k[L]}
\newcommand{\ZL}{\ZZ[L]}

\newcommand{\ab}[1]{{#1}^\text{\rm ab}}
\newcommand{\cm}{Cohen-Ma\-caulay}

\newcommand{\inert}[2]{I_{#1}(#2)}

\newcommand{\tworef}[1]{#1^{(2)}}

\begin{document}

%
%

\title[Cohen-Macaulay rings of invariants]%
{On the Cohen-Macaulay property of multiplicative invariants}

\author{Martin Lorenz}
\address{Department of Mathematics, Temple University,
    Philadelphia, PA 19122}
\email{lorenz@math.temple.edu}
\thanks{Research of the author supported in part by NSF grant DMS-9988756}


\begin{abstract}
We investigate the Cohen-Macaulay property for rings of invariants under
multiplicative actions of a finite group $\G$. By definition, these are
$\G$-actions on Laurent polynomial algebras $\k[x_1^{\pm 1},\dots,x_n^{\pm 1}]$
that stabilize the multiplicative group consisting of all monomials in the
variables $x_i$. For the most part, we concentrate on the case where the base
ring $\k$ is $\ZZ$. Our main result states that if $\G$ acts non-trivially and
the invariant ring
$\ZZ[x_1^{\pm 1},\dots,x_n^{\pm 1}]^\G$ is {\cm} then the abelianized isotropy
groups $\ab{\G}_m$ of all monomials $m$ are generated by the bireflections in
$\G_m$ and at least one $\ab{\G}_m$ is non-trivial.
As an application, we prove the multiplicative version
of Kemper's $3$-copies conjecture.
\end{abstract}

\maketitle



\section*{Introduction}

This article is a sequel to \cite{pathak}.  Unlike in \cite{pathak}, however,
our focus will be specifically on multiplicative invariants.
In detail, let $L \cong \ZZ^n$ denote a lattice on which a finite group $\G$ acts by
automorphisms. The $\G$-action on $L$ extends uniquely to an action by $\k$-algebra
automorphisms on the group algebra $\kL \cong \k[x_1^{\pm 1},\dots,x_n^{\pm 1}]$
over any commutative base ring $\k$. We are interested in the question when the subalgebra
$\kL^\G$ consisting of all $\G$-invariant elements of $\kL$ has the {\cm} property.
The reader is assumed to have some familiarity with {\cm} rings; a good reference on this 
subject is \cite{bruns}.

It is a standard fact that $\kL$ is {\cm} precisely if $\k$ is. On the other hand, while
$\kL^\G$ can only be {\cm} when $\k$ is so, the latter condition is far from sufficient
and rather stringent additional conditions on the action of $\G$ on $L$ are required to
ensure that $\kL^\G$ is {\cm}. Remarkably, the question whether or not
$\kL^\G$ is {\cm}, for any given base ring $\k$, depends only on the rational
isomorphism class of the lattice $L$, that is, the isomorphism class of $L \otimes_{\ZZ}\QQ$ 
as $\QQ[\G]$-module; see Proposition~\ref{SS:rational} below. This is in striking contrast with
most other ring theoretic properties of $\kL^\G$ (e.g., regularity, structure of the class 
group) which tend to be sensitive to the $\ZZ$-type of $L$. For an
overview, see \cite{MIT}.

We will largely concentrate on the case where the base ring $\k$ is $\ZZ$.
This is justified in part by the fact that if $\ZL^\G$ is {\cm} then likewise is
$\kL^\G$ for any {\cm} base ring $\k$ (Lemma~\ref{SS:coefficients}).
Assuming $\ZL^\G$ to be {\cm}, we aim to derive group theoretical consequences for the
isotropy groups $\G_m = \{g \in \G \mid g(m) = m \}$ with $m \in L$.
An element $g \in \G$ will be called a \emph{$k$-reflection} on $L$ if the sublattice
$[g,L] = \{ g(m) - m \mid m \in L \}$ of $L$ has rank at most $k$ or, equivalently,
if the $g$-fixed points of the $\QQ$-space $L \otimes_{\ZZ}\QQ$ have codimension at most $k$. 
As usual, $k$-reflections with $k = 1$ and $k = 2$ will be referred to as \emph{reflections}
and \emph{bireflections}.
For any subgroup $\H \le \G$, we let $\tworef{\H}$ denote the subgroup generated 
by the elements of $\H$ that act as bireflections on $L$.
Our main result now reads as follows.

\begin{thm}
Assume that $\ZL^{\G}$ is {\cm}. Then  $\G_m/\tworef{\G_m}$
is a perfect group (i.e.,  equal to its commutator subgroup) for all $m \in L$.
If $\G$ acts non-trivially on $L$ then some $\G_m$ is non-perfect. 
\end{thm}

It would be interesting to determine if the conclusion of the theorem can be strengthened
to the effect that all isotropy groups $\G_m$ are in fact generated by bireflections on $L$.
I do not know if, for the latter to occur, it is sufficient that $\G$ is  generated by
bireflections. The corresponding fact for reflection groups is known to be true: if $\G$ is
generated by reflections on $L$ (or, equivalently, on $L \otimes_{\ZZ}\QQ$) then
so are all isotropy groups $\G_m$; see
\cite[Theorem 1.5]{rS} or \cite[Exercise 8(a) on p.~139]{nB68}.

There is essentially a complete classification
of finite linear groups generated by bireflections. In arbitrary characteristic, this is due to
Guralnick and Saxl \cite{guralnick}; for the  case of characteristic zero, see
Huffman and Wales \cite{huffman}.
Bireflection groups have been of interest in connection with the problem
of determining all finite linear groups  whose algebra of polynomial invariants
is a complete intersection. Specifically, suppose that $\G \le \GL(V)$ for some finite-dimensional
vector space 
$V$ and let $\cO(V) = \Sym(V^*)$ denote the algebra of polynomial functions on $V$.
It was shown by Kac and Watanabe \cite{vKkW82} and independently by Gordeev \cite{nG82} that 
if the algebra $\cO(V)^\G$ of all $\G$-invariant polynomial functions is
a complete intersection then $\G$ is generated by bireflections on $V$. The classification 
of all groups $\G$ so that $\cO(V)^\G$ is a complete intersection has been
achieved by Gordeev \cite{nG87} and by Nakajima \cite{hN84}.

The last assertion of the above Theorem implies in particular that if $\ZL^\G$ is {\cm} and 
$\G$ acts non-trivially on $L$ then some element of $\G$ acts as a non-trivial bireflection on $L$.
Hence we obtain the following multiplicative version
of Kemper's $3$-copies conjecture:

\begin{cor}
If $\G$ acts non-trivially on $L$ and $r \ge 3$ then $\ZZ[L^{\oplus r}]^{\G}$ is not {\cm}.
\end{cor}

The $3$-copies conjecture was formulated by Kemper \cite[Vermutung 3.12]{gK99} in the context of
polynomial invariants. Using the above notation, the original conjecture states that 
if $1 \neq \G \le \GL(V)$ and the characteristic of the base field of $V$ divides the order of $\G$
(``modular case")  then
the invariant algebra $\cO(V^{\oplus r})^{\G}$ will not be {\cm} for any $r \ge 3$. This
is still open. 
The main factors contributing to our success in the multiplicative case are the following:

\begin{itemize}
\item Multiplicative actions are permutation actions: $\G$ permutes
    the $\k$-basis of $\kL$ consisting of all ``monomials", corresponding to the elements
    of the lattice $L$. Consequently, the cohomology $H^*(\G,\kL)$ is simply the
    direct sum of the various $H^*(\G_m,\k)$ with $m$ running over a transversal for
    the $\G$-orbits in $L$.
\item Up to conjugacy, there are only finitely many finite subgroups of $\GL_n(\ZZ)$ and 
    these groups are explicitly known for small $n$. A crucial observation
    for our purposes is the following:
    if $\G$ is a nontrivial finite perfect subgroup of $\GL_n(\ZZ)$ such that no $1 \neq g \in \G$ 
    has eigenvalue $1$ then $\G$ is isomorphic to the binary icosahedral group and $n \ge 8$; see 
    Lemma~\ref{SS:fpf} below.
\end{itemize}

A brief outline of the contents of the this article is as follows. The short preliminary 
Section~\ref{S:general} is devoted to general actions of a finite group $\G$ on a commutative 
ring $R$. This material relies rather heavily on \cite{pathak}. We liberate a technical result from 
\cite{pathak} from any a priori hypotheses on the characteristic; the new version (Proposition~\ref{SS:necessary})
states that if $R$ and $R^\G$ are both {\cm} and $H^i(\G,R) = 0$ for $0 < i < k$ then
$H^k(\G,R)$ is detected by $k+1$-reflections.
Section~\ref{S:mult} then specializes to the case of multiplicative actions. We assemble the main
tools required for the proof of the Theorem, which is presented in Section~\ref{S:cm}.
The article concludes with a brief discussion of possible avenues for further investigation and 
some examples. 




\section{Finite Group Actions on Rings} \label{S:general}


\subsection{} \label{SS:generalintro}
In this section, $R$ will be a commutative
ring on which a finite group $\G$ acts by ring automorphisms
$r \mapsto g(r)$ $(r \in R, g \in \G)$. The subring
of $\G$-invariant elements of $R$ will be denoted by $R^\G$. 


\subsection{Generalized reflections} \label{SS:reflections}
Following \cite{gordeev}, we will say an element $g \in \G$ acts as a 
\emph{$k$-reflection} on $R$ if $g$ belongs to the inertia group
$$
\inert{\G}{\fP} = \{ g \in \G \mid g(r) - r \in \fP\ \forall r \in R\}
$$
of some prime ideal $\fP \in \Spec R$ with $\height \fP \le k$.
The cases $k = 1$ and $k = 2$ will be referred to as \emph{reflections}
and \emph{bireflections}, respectively. Define the ideal $I_R(g)$ of $R$ by
\begin{equation*} \label{E:ig}
I_R(g) = \sum_{r \in R} (g(r) - r)R \ .
\end{equation*}
Evidently, $\fP \supseteq I_R(g)$
is equivalent to $g \in \inert{\G}{\fP}$. Thus:
$$
\text{$g$ is a $k$-reflection on $R$
if and only if $\height I_R(g) \le k$.}
$$
For each subgroup $\H \le \G$, we put
\begin{equation*} \label{E:iH}
I_R(\H) = \sum_{g \in \H} I_R(g) \ .
\end{equation*}
It suffices to  let $g$ run over a set of generators of the group $\H$ in this sum.


\subsection{A height estimate} \label{SS:height}
The cohomology $H^*(\G,R) = \oplus_{n \ge 0} H^n(\G,R)$ has a canonical
$R^\G$-module structure:
for each $r\in R^\G$, the map $\rho \colon R \to R$, $s \mapsto rs$, is
$\G$-equivariant and hence it induces a map on cohomology $\rho_*: H^*(\G,R)
\to H^*(\G,R)$. The element $r$ acts on $H^*(\G,R)$ via $\rho_*$. Let
$\res^\G_\H \colon H^*(\G,R) \to H^*(\H,R)$ denote the restriction map.

The following lemma extends \cite[Proposition 1.4]{pathak}.

\begin{lem} \label{L:height}
For any $x\in H^*(\G,R)$,
$$
\height\ann_{R^\G}(x) \ge \inf\{\height I_R(\H) \mid 
\H \le \G, \res^\G_\H(x) \neq 0 \}\ .
$$
\end{lem}

\begin{proof}
Put $\cX=\{ \H \le \G \mid \res^\G_\H(x) = 0 \}$. For each $\H \le \G$,
let $R_\H^\G$ denote the image of the relative trace map 
$R^\H \to R^\G$, $r \mapsto \sum_{g} g(r)$,
where $g$ runs over a transversal for the cosets $g\H$ of $\H$ in $\G$.
By \cite[Lemma 1.3]{pathak},
$$
R^\G_\H \subseteq \ann_{R^\G}(x) \quad \text{for all $\H \in \cX$.}
$$
To prove the lemma, we may assume that $\ann_{R^\G}(x)$ is a proper ideal of
$R^\G$; for, otherwise $\height\ann_{R^\G}(x) = \infty$. Choose a prime ideal 
$\fp$ of $R^\G$ with $\fp \supseteq \ann_{R^\G}(x)$ and $\height \fp = 
\height \ann_{R^\G}(x)$. If $\fP$ is a prime of $R$ that lies over $\fp$ then 
$$
R^\G_\H \subseteq \fP \quad \text{for all $\H \in \cX$}
$$
and $\height \fP =  \height \fp$. By \cite[Lemma 1.1]{pathak},
the above inclusion implies that
$$
[\inert{\G}{\fP} : \inert{\H}{\fP}] \in \fP \quad \text{for all $\H \in \cX$.}
$$
Put $p = \ch R/\fP$ and let $\cP \le \inert{\G}{\fP}$ be a Sylow $p$-subgroup of
$\inert{\G}{\fP}$ (so $\cP = 1$ if $p = 0$). Then $I_R(\cP) \subseteq \fP$ and 
$[\inert{\G}{\fP} : \cP] \notin \fP$.
Hence, $\cP \notin \cX$ and $\height I_R(\cP) \le 
\height \fP = \height \ann_{R^\G}(x)$. This proves the lemma.
\end{proof}

We remark that the lemma and its proof carry over verbatim to the more
general situation where $H^*(\G,R)$ is replaced by $H^*(\G,M)$, where $M$ is
some module over the skew group ring of $\G$ over $R$; cf.~\cite{pathak}.
However, we will not be concerned with this generalization here.


\subsection{A necessary condition} \label{SS:necessary}

In this section, we assume  
that $R$ is noetherian as $R^\G$-module. This assumption is satisfied
whenever $R$ is an affine algebra over some noetherian subring 
$\k \subseteq R^\G$; see \cite[Th\'eor\`eme 2 on p.~33]{nB64}. 
Put 
\begin{equation} \label{E:Xk}
\cX_k = \{ \H \le \G \mid  \height I_R(\H) \le k \} \ .
\end{equation}
Note that each
$\H \in \cX_k$ consists of $k$-reflections on $R$. The following proposition is
a characteristic-free version of \cite[Proposition 4.1]{pathak}.

\begin{prop} 
Assume that $R$ and $R^{\G}$ are {\cm}. If $H^i(\G,R) = 0$ $(0 < i < k)$ then
the restriction map
$$
\res^{\G}_{\cX_{k+1}} \colon H^k(\G,R) \to \prod_{\H \in \cX_{k+1}}
H^k(\H,R)
$$
is injective.
\end{prop}

\begin{proof}
We may assume that $H^k(\G,R) \neq 0$. Let $x \in H^k(\G,R)$ be nonzero and put
$\fa = \ann_{R^\G}(x)$. By \cite[Proposition 3.3]{pathak}, $\depth \fa \le k+1$.
Since $R^\G$ is {\cm}, $\depth \fa = \height \fa$. Thus, Lemma~\ref{L:height} implies that
$k+1 \ge \height I_R(\H)$ for some $\H \le \G$ with $\res^\G_\H(x) \neq 0$.
The proposition follows.
\end{proof}

Note that the vanishing hypothesis on $H^i(\G,R)$ is vacuous for $k=1$. Thus, $H^1(\G,R)$ is
detected by bireflections whenever $R$ and $R^{\G}$ are both {\cm}.


\section{Multiplicative Actions} \label{S:mult}


\subsection{} \label{SS:multintro}

For the remainder of this article, $L$ will denote a lattice on which the finite 
group $\G$ acts by automorphisms $m \mapsto g(m)$ $(m \in L, g \in \G)$. The
group algebra of $L$ over some commutative base ring $\k$ will be denoted by $\kL$. 
We will use additive notation in $L$. The $\k$-basis element of
$\kL$ corresponding to the lattice element $m \in L$ will be written as
$$
\x^m \ ;
$$
so $\x^0 = 1$, $\x^{m+m'} = \x^m \x^{m'}$, and $\x^{-m} = (\x^m)^{-1}$. 
The action of $\G$ on $L$ extends uniquely to an action by $\k$-algebra
automorphisms on $\kL$:
\begin{equation*} \label{E:multact}
g(\sum_{m \in L} k_m\x^m) = \sum_{m \in L} k_m\x^{g(m)} \ .
\end{equation*}
The invariant algebra $\kL^\G$ is a free $\k$-module: a $\k$-basis  is
given by the $\G$-orbit sums $\osum(m) = \sum_{m' \in \G(m)} \x^{m'}$, 
where $\G(m)$ denotes the $\G$-orbit of $m \in L$. Since all orbit
sums are defined over $\ZZ$, we have
\begin{equation} \label{E:overZ}
\kL^{\G} = \k \otimes_{\ZZ}\ZL^{\G}\ .
\end{equation}


\subsection{} \label{SS:I(H)}
Let $\H$ be a subgroup of $\G$. We compute the height of the ideal 
$I_{\kL}(\H)$ from \S\ref{SS:reflections}.
Let
$$
L^{\H} = \{ m \in L \mid g(m) = m \text{ for all $g \in \H$}\}
$$
denote the lattice of $\H$-invariants in $L$ and
define the sublattice $[\H,L]$ of $L$ by
$$
[\H,L] = \sum_{g \in \H}\, [g,L]\ ,
$$
where $[g,L] = \{ g(m) - m \mid m \in L \}$.
It suffices to  let $g$ run over a set of generators of the 
group $\H$ in the above formulas.

\begin{lem} 
With the above notation, $\kL/I_{\kL}(\H) \cong \k[L/[\H,L]]$ and
$$
\height I_{\kL}(\H) = \rank [\H,L] = \rank L - \rank L^\H \ .
$$
\end{lem}

\begin{proof}
Since the ideal $I_{\kL}(\H)$ is generated
by the elements $\x^{g(m) - m} - 1$ with $m \in L$ and $g \in \H$,
the isomorphism $\kL/I_{\kL}(\H) \cong \k[L/[\H,L]]$ is clear.

To prove the equality $\rank [\H,L] = \rank L - \rank L^\H$, note that
the rational group algebra $\QQ[\H]$ is the direct sum of the ideals 
$\QQ\left(\sum_{g \in \H} g\right)$ and $\sum_{g \in \H} \QQ(g -
1)$. This implies $L \otimes_{\ZZ}\QQ = \left( L^{\H}
\otimes_{\ZZ}\QQ\right) \oplus \left( [\H,L]
\otimes_{\ZZ}\QQ\right)$. Hence, $\rank L = \rank L^\H + \rank
[\H,L]$.

To complete the proof, it suffices to show that
$$
\height \fP = \rank [\H,L] 
$$
holds for any minimal covering prime $\fP$ of $I_{\kL}(\H)$.
Put $A = L/[\H,L]$ and $\bar{\fP} = \fP/I_{\kL}(\H)$, a minimal prime of
$\kL/I_{\kL}(\H) = \k[A]$. Further, put $\fp = \bar{\fP} \cap \k = \fP \cap
\k$. Since the extension $\k \into \k[A] = \kL/I_{\kL}(\H)$ is free, $\fp$
is a minimal prime of $\k$; see \cite[Cor.~on p.~AC
VIII.15]{nB83}. Hence, descending chains of primes in $\kL$
starting with $\fP$ correspond in a $1$-to-$1$ fashion to
descending chains of primes of $Q(\k/\fp)[L]$ starting with the
prime that is generated by $\fP$. Thus, replacing $\k$ by
$Q(\k/\fp)$, we may assume that $\k$ is a field. But then
$$
\height \fP = \dim \kL - \dim \kL/\fP  = \rank L - \dim \kL/\fP \
.
$$
Let $\bar{\fP}_0 = \bar{\fP} \cap \k[A_0]$, where $A_0$ denotes
the torsion subgroup of $A$. Since $\bar{\fP}$ is minimal, we have 
$\bar{\fP} = \bar{\fP}_0\k[A]$ and so
$\kL/\fP \cong \k_0[A/A_0]$, where $\k_0 = \k[A_0]/\bar{\fP}_0$ is
a field. Thus, $\dim \kL/\fP = \rank A/A_0$. Finally, $\rank A/A_0
= \rank A = \rank L - \rank [\H,L]$, which completes the proof.
\end{proof}

Specializing the lemma to the case where $\H = \gen{g}$ for some $g \in \G$,
we see that $g$ acts as a $k$-reflection on $\kL$ if and only if
$g$ acts as a $k$-reflection on $L$, that is,
$$
\rank [g,L] \le k \ .
$$
Moreover, the collection of subgroups $\cX_k$ in equation \eqref{E:Xk}
can now be written as
\begin{equation} \label{E:Xk'}
\cX_k = \{ \H \le \G \mid \rank L/L^\H \le k \} \ .
\end{equation}


\subsection{Fixed-point-free lattices for perfect groups} \label{SS:fpf}

The $\G$-action on $L$ is called \emph{fixed-point-free} if $g(m) \neq m$ holds
for all $0 \neq m \in L$ and $1\neq g \in \G$. Recall also that the group
$\G$ is said to be \emph{perfect} if $\ab{\G} = \G/[\G,\G] = 1$.

\begin{lem}
Assume that $\G$ is a nontrivial perfect group acting fixed-point-freely on the 
nonzero lattice $L$.
Then $\G$ is isomorphic to the binary icosahedral group $2.\A_5 \cong \SL_2(\FF_5)$ 
and $\rank L$ is a multiple of $8$.
\end{lem}

\begin{proof}
Put $V = L \otimes_{\ZZ} {\CC}$, a nonzero fixed-point-free $\CC[\G]$-module.
By a well-known theorem of Zassenhaus (see \cite[Theorem 6.2.1]{wolf}), $\G$ is
isomorphic to the binary icosahedral group $2.\A_5$ and the irreducible constituents of 
$V$ are $2$-dimensional. The binary icosahedral group has two irreducible
complex representations of degree $2$; they are Galois conjugates of each other and
both have Frobenius-Schur indicator $-1$.
We denote the corresponding $\CC[\G]$-modules by $V_1$ and $V_2$. Both $V_i$ occur with 
the same multiplicity in $V$, since $V$ is defined over $\QQ$. Thus,
$V \cong \left( V_1 \oplus V_2 \right)^m$ for some $m$ and $\rank L = 4m$. We have to
show that $m$ is even. Since both $V_i$ have indicator $-1$,
it follows that $V_1 \oplus V_2$ is not defined over $\RR$, whereas each $V_i^2$ is defined over $\RR$;
see \cite[(9.21)]{isaacs}. Thus, $V_1 \oplus V_2$ represents an element $x$ of order $2$ in the
cokernel of the scalar extension map $G_0(\RR[\G]) \to G_0(\CC[\G])$, and $mx = 0$.
Therefore, $m$ must be even, as desired.
\end{proof}

We remark that the binary icosahedral group 
$2.\A_5$ is isomorphic to the subgroup of the nonzero quaternions $\HH^*$
that is generated by  $(a + i + j a^*)/2$ and $(a + j + k a^*)/2$, where 
$a = (1 + \sqrt{5})/2$ and $a^* = (1 - \sqrt{5})/2$ and $\{1,i,j,k\}$ is the standard $\RR$-basis
of $\HH$. Thus, letting $2.\A_5$ act on 
$\HH$ via left multiplication, $\HH$ becomes a $2$-dimensional fixed-point-free complex representation
of $2.\A_5$. It is easy to see that this representation can be
realized over $K = \QQ(i,\sqrt{5})$; so $\HH = V \otimes_K \CC$ with $\dim_{\QQ}V = 2[K:\QQ] = 8$.
Any $2.\A_5$-lattice for $V$ will be fixed-point-free and have rank $8$.


\subsection{Isotropy groups} \label{SS:isotropy}

The isotropy group of an element $m \in L$ in $\G$ will be denoted by $\G_m$; so
$$
\G_m = \{ g \in \G \mid g(m) = m \} \ .
$$
The $\G$-lattice $L$ is called \emph{faithful} if $\Ker_{\G}(L) = \bigcap_{m \in L} \G_m = 1$.
The following lemma, at least part (a), is well-known. We include the proof for the
reader's convenience.

\begin{lem}
\begin{enumerate}
\item The set of isotropy groups $\{ \G_m \mid m \in L \}$ is closed under
conjugation and under taking intersections.
\item Assume that the $\G$-lattice $L$ is faithful. If $\G_m \ (m \in L)$ is a  minimal non-identity
isotropy group  then
$\G_m$ acts fixed-point-freely on $L/L^{\G_m}\neq 0$.
\end{enumerate}
\end{lem}

\begin{proof}
Consider the $\QQ[\G]$-module $V = L \otimes_{\ZZ} {\QQ}$. The collection
of isotropy groups $\G_m$ remains unchanged when allowing $m \in V$. Moreover, for any
subgroup $\H \le \G$,  $L/L^{\H}$ is an $\H$-lattice with $L/L^{\H} \otimes_{\ZZ}\QQ \cong V/V^{\H}$. 

(a) The first assertion is clear, since ${^g\G_m} = \G_{g(m)}$ holds for all $g \in \G, m \in V$.
For the second assertion, let $M$ be a non-empty subset of $V$ and put $\G_M = \bigcap_{m \in M} \G_m$.
We must show that $\G_M = \G_m$ for some $m \in V$. 
Put $W = V^{\G_M}$. 
If $g \in \G \setminus \G_M$ then $W^g = \{ w \in W \mid g(w) = w \}$ is a proper subspace of $W$, since some 
element of $M$ does not belong to $W^g$. 
Any $m \in W \setminus \bigcup_{g \in \G \setminus \G_M} W^g$ satisfies $\G_m = \G_M$.

(b) Let $\H = \G_m$ be a minimal non-identity member of $\{ \G_m \mid m \in V \}$. As $\QQ[\H]$-modules,
we may identify $V$ and $V^{\H} \oplus V/V^{\H}$ . If $0 \neq v \in V/V^{\H}$
then $\H_v = \H \cap \G_v \subsetneq \H$. In view of (a), our minimality assumption on $\H$ forces
$\H_v = 1$. Thus, $\H$ acts fixed-point-freely on $V/V^\H$, and hence on $L/L^\H$.
\end{proof}

\begin{prop}
Assume that $L$ is a faithful $\G$-lattice such that all minimal isotropy groups $1 \neq \G_m \ (m \in L)$
are perfect. Then $\rank L/L^{\H} \ge 8$ holds for every nonidentity subgroup $\H \le \G$.
\end{prop}

\noindent
In the notation of equation \eqref{E:Xk'}, the conclusion of the proposition can be stated as follows:
$$
\cX_k = \{ 1\}\ \text{for all $k < 8$.}
$$

\begin{proof}[Proof of the Proposition]
Put $\bar{\H} = \bigcap_{m \in L^\H} \G_m$. Then $\bar{H} \supseteq \H$ and
$L^{\bar{\H}} = L^\H$. Lemma~\ref{SS:isotropy}(a) further implies that $\bar{\H} = \G_m$
for some $m$. Replacing $\H$ by $\bar{H}$, we may assume that $\H$ is a nonidentity isotropy group. 
If $\H$ is not minimal then replace $\H$
by a smaller nonidentity isotropy group; this does not increase the value of $\rank L/L^{\H}$.
Thus, we may assume that $\H$ is a minimal nonidentity isotropy group, and hence $\H$ is perfect. By
Lemma~\ref{SS:isotropy}(b), $\H$ acts fixed-point-freely on $L/L^\H \neq 0$
and Lemma~\ref{SS:fpf} implies that $\rank L/L^{\H} \ge 8$, proving the proposition.
\end{proof}


\subsection{Cohomology} \label{SS:restriction}
Let $\cX$ denote any collection of subgroups of $\G$ that is closed
under conjugation and under taking subgroups. 
We will investigate injectivity of the restriction map
$$
\res^{\G}_{\cX} \colon H^k(\G,\kL) \to \prod_{\H \in \cX} H^k(\H,\kL) \ .
$$
This map was considered in  Proposition~\ref{SS:necessary} for $\cX = \cX_{k+1}$.

\begin{lem} 
The map $\res^{\G}_{\cX}  \colon H^k(\G,\kL) \to \prod_{\H \in \cX} H^k(\H,\kL)$ is
injective if and only if the restriction maps
$$
H^k(\G_m,\k) \to \prod_{\substack{\H \in \cX \\ \H \le \G_m}} H^k(\H,\k)
$$
are injective for all $m \in L$.
\end{lem}

\begin{proof} As $\k[\G]$-module, $\kL$ is a permutation module:
$$
\kL \cong \bigoplus_{m \in \G\backslash L} \k[\G/\G_m] \ ,
$$
where $\k[\G/\G_m] = \k[\G]\otimes_{\k[\G_m]}\k$ and $\G\backslash L$ is a transversal for the 
$\G$-orbits in $L$. For each subgroup $\H \le \G$,
$$
\k[\G/\G_m]\big|_{\H} \cong \bigoplus_{g \in \H\backslash \G/\G_m} \k[\H/{^g\G_m}\cap \H]\ ;
$$
see \cite[10.13]{curtis}. Therefore, $\res^{\G}_{\H}$ is the direct sum
of the restriction maps 
$$
H^k(\G,\k[\G/\G_m]) \to  H^k(\H,\k[\G/\G_m]) =
\bigoplus_{g \in \H\backslash \G/\G_m} H^k(\H,\k[\H/{^g\G_m}\cap \H])\ .
$$
By the Eckmann-Shapiro Lemma \cite[III(5.2),(6.2)]{brown},
$H^k(\G,\k[\G/\G_m]) \cong  H^k(\G_m,\k)$ and $H^k(\H,\k[\H/{^g\G_m}\cap \H])
\cong H^k({^g\G_m}\cap \H,\k)$.
In terms of these isomorphisms, the above restriction map becomes
\begin{alignat*}{3}
\rho_{\H,m} \colon H^k&(\G_m,\k) & \quad &\to & \quad
\bigoplus_{g \in \H\backslash \G/\G_m} H^k&({^g\G_m}\cap \H,\k) \\
&[f] & \quad &\mapsto & \quad ([ \underline{h} \mapsto &f(g^{-1}\underline{h}g) ])_g
\end{alignat*}
where $[\,.\,]$ denotes the cohomology class of a $k$-cocycle and $\underline{h}$
stands for a $k$-tuple of elements of ${^g\G_m}\cap \H$.
Therefore,
$$
\Ker \rho_{\H,m} = \bigcap_{g \in \H\backslash \G/\G_m} \Ker \left( \res^{\G_m}_{\G_m\cap\H^g}
\colon H^k(\G_m,\k) \to H^k(\G_m\cap\H^g,\k) \right)\ .
$$
Thus, $\Ker\res^{\G}_{\cX}$ is isomorphic to the direct sum of the kernels of the
restriction maps
$$
H^k(\G_m,\k) \to \prod_{\H \in \cX} H^k(\G_m\cap\H^g,\k)
$$
with $m \in \G\backslash L$. Finally, by hypothesis on $\cX$, the groups $\G_m\cap\H^g$ with
$\H \in \cX$ are exactly the
groups $\H \in \cX$ with $\H \le \G_m$. The lemma follows.
\end{proof}

\begin{cor}
Let $\k = \ZZ/(|\G|)$ and $k = 1$. Then
$\res^{\G}_{\cX}$ 
injective if and only if all $\ab{\G}_m$ $(m \in L)$ are generated by the images of the subgroups 
$\H \le \G_m$ with $\H \in \cX$.
\end{cor}

\begin{proof}
By the lemma with $k=1$, the hypothesis on the restriction map says that all  restrictions
$$
H^1(\G_m,\k) \to \prod_{\substack{\H \in \cX \\ \H \le \G_m}} H^1(\H,\k)
$$
are injective. Now, for each $\H \le \G$, $H^1(\H,\k) = \Hom(\ab{\H},\k) \cong \ab{\H}$, where
the last isomorphism holds by our choice of $\k$. Therefore, injectivity of the above map
is equivalent to $\ab{\G}_m$ being generated by the images of all $\H \le \G_m$ with $\H \in \cX$.
\end{proof}


\section{The {\cm} Property} \label{S:cm}


\subsection{} \label{SS:cmintro}
Continuing with the notation of \S\ref{SS:multintro}, we now turn to the question
when the invariant algebra $\kL^\G$ is {\cm}. 
Our principal tool will be Proposition~\ref{SS:necessary}.
We remark that the {\cm} hypothesis of Proposition~\ref{SS:necessary}
simplifies slightly in the setting of multiplicative actions: it suffices to assume that
$\kL^\G$ is {\cm}. Indeed, in this case the base ring $\k$ is also {\cm},
because $\kL^\G$ is free over $\k$, and then $\kL$ is {\cm} as well;
see \cite[Exercise 2.1.23 and Theorems 2.1.9, 2.1.3(b)]{bruns}.


\subsection{Base rings} \label{SS:coefficients}

Our main interest is in the case where $\k = \ZZ$. As the following lemma shows,
if $\ZL^\G$ is {\cm} then so is $\kL^\G$ for any {\cm} base ring $\k$.

\begin{lem} 
The following are equivalent:
\begin{enumerate}
\item $\ZL^\G$ is {\cm};
\item $\kL^\G$ is {\cm} whenever $\k$ is;
\item $\kL^\G$ is {\cm} for $\k = \ZZ/(|\G|)$;
\item $\FF_p[L]^\G$ is {\cm} for all primes $p$ dividing $|\G|$.
\end{enumerate}
\end{lem}

\begin{proof}
(a) $\Rightarrow$ (b): Put $S = \kL^\G$ and consider the extension of
rings $\k \into S$. This extension is free; see \S\ref{SS:multintro}. 
By \cite[Exercise 2.1.23]{bruns},
$S$ is {\cm} if (and only if) $\k$ is {\cm} and, for all $\fP \in \Spec S$, the
fibre $S_{\fP}/{\fp}S_{\fP}$ is {\cm}, where $\fp = \fP \cap \k$. But
$S_{\fP}/{\fp}S_{\fP}$ is a localization of $(S/{\fp}S)_{\fp \setminus 0} \cong
Q(\k/\fp)[L]^\G$; see equation \eqref{E:overZ}. Therefore, by 
\cite[Theorem 2.1.3(b)]{bruns}, it suffices to show that $Q(\k/\fp)[L]^\G$ is
{\cm}. In other words, we may assume that $\k$ is a field. By 
\cite[Theorem 2.1.10]{bruns}, we may further assume that $\k = \QQ$ or $\k = \FF_p$.
But equation \eqref{E:overZ} implies that $\QQ[L]^\G = \ZL^\G_{\ZZ \setminus 0}$ 
and $\FF_p[L]^\G \cong \ZL^{\G}/(p)$. Since $\ZL^\G$ is assumed {\cm}, 
\cite[Theorem 2.1.3]{bruns} implies that $\QQ[L]^\G$ and $\FF_p[L]^\G$ are {\cm}, 
as desired.

(b) $\Rightarrow$ (c) is clear.

(c) $\Rightarrow$ (d): Write $|\G| = \prod_p p^{n_p}$. Then 
$\kL \cong \prod_p \ZZ/(p^{n_p})[L]^\G$
and $\ZZ/(p^{n_p})[L]^\G$ is a localization of $\kL^\G$.
Therefore, $\ZZ/(p^{n_p})[L]^\G$ is {\cm}, by \cite[Theorem 2.1.3(b)]{bruns}.
If $n_p \neq 0$ then it follows from \cite[Theorem 2.1.3(a)]{bruns} that
$\ZZ_{(p)}[L]^\G$ and $\FF_p[L]^\G \cong \ZZ_{(p)}[L]^\G/(p)$ are {\cm}.

(d) $\Rightarrow$ (a): First, (d) implies that $\FF_p[L]^\G$ is {\cm} for all primes $p$.
For, if $p$ does not divide $|\G|$ then
$\FF_p[L]^\G$ is always {\cm}; see \cite[Corollary 6.4.6]{bruns}. Now let $\fP$
be a maximal ideal of $\ZL$. Then $\fP \cap \ZZ = (p)$ for some prime $p$ and
$\ZL^{\G}_{\fP}/(p)$ is a localization of $\ZL^{\G}/(p) = \FF_p[L]^\G$.
Thus, $\ZL^{\G}_{\fP}/(p)$ is {\cm} and \cite[Theorem 2.1.3(a)]{bruns} further
implies that $\ZL^{\G}_{\fP}$ is {\cm}. Since, $\fP$ was arbitrary, (a) is proved.
\end{proof}

Since normal rings of (Krull) dimension at most $2$ are {\cm}, the implication (d)
$\Rightarrow$ (b) of the lemma shows that $\kL^\G$ is certainly {\cm} whenever
$\k$ is {\cm} and $L$ has rank at most $2$.


\subsection{Proof of the Theorem} \label{SS:theorem}

We are now ready to prove the Theorem stated in the Introduction.
Recall that, for any subgroup $\H \le \G$, $\tworef{\H}$ denotes the subgroup generated
by the elements of $\H$
that act as bireflections on $L$ or, equivalently, by the subgroups of $\H$ that belong to
$\cX_2$; see \eqref{E:Xk'}. Throughout, we assume that $\ZL^{\G}$ is {\cm}.

We first show that  $\G_m/\tworef{\G_m}$ 
is a perfect group for all $m \in L$. Put $\k = \ZZ/(|\G|)$. Then
$\kL^{\G}$ is {\cm}, by Lemma~\ref{SS:coefficients}.
Therefore, the restriction $H^1(\G,\kL) \to \prod_{\H \in \cX_2} H^1(\H,\kL)$ is 
injective, by Proposition~\ref{SS:necessary}; see the remark in \S\ref{SS:cmintro}.
Corollary~\ref{SS:restriction} yields that all $\ab{\G}_m$ 
are generated by the images of the subgroups $\H \le \G_m$ with $\H \in \cX_2$.
In other words, each $\ab{\G}_m$ is generated by the images of the bireflections in $\G_m$.
Therefore, $\ab{\left(\G_m/\tworef{\G_m}\right)} = 1$, as desired.

Now assume that $\G$ acts non-trivially on $L$. Our goal is to show that
some isotropy group $\G_m$ is non-perfect. Suppose otherwise.
Replacing $\G$ by  $\G/\Ker_{\G}(L)$ we may assume that 
$1 \neq \G$ acts faithfully on $L$. 
Then $\cX_k = \{1\}$ for all $k < 8$, by Proposition~\ref{SS:isotropy}.
It follows that 
$$
k = \inf\{ i > 0 \mid H^i(\G, \kL) \neq 0 \} \ge 7 \ .
$$
Indeed, if $k < 7$ then Proposition~\ref{SS:necessary} implies that
$0 \neq H^k(\G,\kL)$ embeds into $\prod_{\H \in \cX_{k+1}}
H^k(\H,\kL)$ which is trivial, because $\cX_{k+1} = \{1\}$.
By Lemma~\ref{SS:restriction} with $\cX = \{ 1 \}$, we conclude that 
$$
H^i(\G_m, \k) = 0 \ \text{for all $m \in L$ and all $0 < i < 7$.}
$$
On the other hand, choosing $\G_m$ minimal with 
$\G_m \neq 1$, we know by Lemmas~\ref{SS:fpf} and \ref{SS:isotropy}(b) that $\G_m$ is isomorphic to the 
binary icosahedral group $2.\A_5$. 
The cohomology of $2.\A_5$ is $4$-periodic (see \cite[p.~155]{brown}). Hence,
$H^3(\G_m,\k) \cong H^{-1}(\G_m,\k) = \ann_{\k}(\sum_{\G_m}g) \cong \ZZ/(|\G_m|) \neq 0$.
This contradiction completes the proof of the Theorem.
\qed





\subsection{Rational invariance}
\label{SS:rational}

We now show that the {\cm} property of $\kL^\G$ depends only on the rational
isomorphism class of the $\G$-lattice $L$. Recall that $\G$-lattices $L$ and $L'$ are said to be
\emph{rationally isomorphic} if $L \otimes_{\ZZ}\QQ \cong L'\otimes_{\ZZ}\QQ$ as $\QQ[\G]$-modules.
In this section, $\k$ denotes any commutative base ring.

\begin{prop} \label{P:rational}
If $\kL^{\G}$ is {\cm} then so is  $\k[L']^{\G}$ for any $\G$-lattice $L'$ that is rationally 
isomorphic to $L$.
\end{prop}

\begin{proof} Assume that $L \otimes_{\ZZ}\QQ \cong L'\otimes_{\ZZ}\QQ$. Replacing $L'$ by an
isomorphic copy inside $L \otimes_{\ZZ}\QQ$, we may assume that $L \supseteq L'$ and $L/L'$ is
finite. Then $\kL$ is finite over $\k[L']$ which in turn is integral over
$\k[L']^\G$. Therefore, $\kL$ is integral over $\k[L']^\G$, and hence so is $\kL^{\G}$.

We now invoke a ring-theoretic result of Hochster and Eagon \cite{mHjE71} 
(or see \cite[Theorem 6.4.5]{bruns}):
Let $R \supseteq S$ be an integral extension of commutative rings having a Reynolds operator, that is,
an $S$-linear map $R \to S$ that restricts to the identity on $S$.
If $R$ is {\cm} then $S$ is {\cm} as well.

To construct the requisite Reynolds operator, consider the truncation map
$$
\pi \colon \kL \to \k[L']\ ,\quad \sum_{m \in L} k_m\x^m \mapsto \sum_{m \in L'} k_m\x^m \ .
$$
This is a Reynolds operator for the extension $\kL \supseteq \k[L']$ that satisfies $\pi(g(f)) = g(\pi(f))$
for all $g \in \G$, $f \in \kL$. Therefore, $\pi$ restricts to a Reynolds operator $\kL^{\G} 
\to \k[L']^{\G}$ and the proposition follows.
\end{proof}

The proposition in particular allows to
reduce the general case of the {\cm} problem
for multiplicative invariants to the case of  effective $\G$-lattices. 
Recall that the $\G$-lattice $L$ is  \emph{effective} if $L^{\G} = 0$. For any
$\G$-lattice $L$, the quotient $L/L^{\G}$ is an effective $\G$-lattice;
this follows, for example, from the fact that $L$ is rationally isomorphic to 
the $\G$-lattice $L^\G \oplus L/L^\G$.

\begin{cor} \label{C:rational}
$\kL^{\G}$ is {\cm} if and only if this holds for $\k[L/L^{\G}]^{\G}$.
\end{cor}

\begin{proof}
By the proposition, we may replace $L$ by $L' = L^\G \oplus L/L^\G$. But
$\k[L']^G \cong  \k[L/L^\G]^\G \otimes_{\k} \k[L^\G]$, a Laurent polynomial algebra 
over $k[L/L^\G]^\G$.
Thus, by \cite[Theorems 2.1.3 and 2.1.9]{bruns}, 
$\k[L']^\G$ is {\cm} if and only if $\k[L/L^\G]^\G$ is {\cm}.
The corollary follows.
\end{proof}


\subsection{Remarks and examples}
\label{SS:concluding}

\subsubsection{Abelian bireflection groups} \label{SSS:abelian}
It is not hard to show that if $\G$ is a finite abelian group acting as
a bireflection group on the lattice $L$ then $\ZL^\G$ is {\cm}. Using
Corollary~\ref{SS:rational} and an induction on $\rank L$, the proof reduces to
the verification that  $\ZL^\G$ is {\cm} for $L = \ZZ^n$  and
$\G = \diag(\pm 1, \dots, \pm 1) \cap \SL_n(\ZZ)$. Direct computation shows that,
for $n \ge 2$, 
$$
\ZL^{\G} = \ZZ[\xi_1,\dots,\xi_n] \oplus \eta \ZZ[\xi_1,\dots,\xi_n]
$$
where $\xi_i = \x^{e_i} + \x^{-e_i}$ is the $\G$-orbit sum of the standard basis
element $e_i \in \ZZ^n$ and $\eta$ is the orbit sum of $\sum_i e_i = (1,\dots,1)$.

It would be worthwhile to try and extend this result to larger classes of bireflection
groups. The aforementioned classification of bireflection groups in 
\cite{guralnick} will presumably be helpful in this endeavor.

\subsubsection{Subgroups of reflection groups} \label{SSS:subgroups}
Assume that $\G$ acts as a reflection group on the lattice $L$ and let $\H$ be a subgroup
of $\G$ with $[\G : \H] =2$. Then $\H$ acts as a bireflection group. (More generally, if
$\G$ acts as a $k$-reflection group and $[\G : \H] = m$ then $\H$ acts as a 
$km$-reflection group; see \cite{MIT}.) Presumably $\ZL^\H$ will always be
{\cm}, but I have no proof. For an explicit example, let $\G = \Sy_n$ be the 
symmetric group on $\{1,\dots,n\}$
and let $L = U_n$ be the standard permutation lattice for $\Sy_n$; so
$U_{n} = \bigoplus_{i=1} ^n\ZZ e_i$ with $s(e_i) = e_{s(i)}$ for $s \in \Sy_n$.
Transpositions act as reflections on $U_n$ and $3$-cycles as bireflections. 
Let $\A_n \le \Sy_n$ denote the alternating group. To compute $\ZZ[U_n]^{\A_n}$,
put $x_i = \x^{e_i} \in \ZZ[U_n]$. Then $\ZZ[U_n] = \ZZ[x_1,\dots,x_n][s_n^{-1}]$, 
where $s_n = \x^{\sum_1^n e_i} = \prod_1^n x_i$
is the $n^\text{th}$ elementary symmetric function,
and $\Sy_n$ acts via $s(x_i) = x_{s(i)}\ (s \in \Sy_n)$. Therefore,
$\ZZ[U_n]^{\A_n} = \ZZ[x_1,\dots,x_n]^{\A_n}[s_n^{-1}]$. The ring
$\ZZ[x_1,\dots,x_n]^{\A_n}$ of polynomial $\A_n$-invariants has the following form; see
\cite[Theorem 1.3.5]{smith}:
$\ZZ[x_1,\dots,x_n]^{\A_n} = \ZZ[s_1,\dots,s_n] \oplus d \ZZ[s_1,\dots,s_n]$,
where $s_i$ is the $i^\text{th}$ elementary symmetric function and
$$
d = \tfrac{1}{2}\left( \Delta + \Delta_+ \right)
$$
with $\Delta_+ = \prod_{i < j} (x_i + x_j)$ and $\Delta = \prod_{i < j} (x_i - x_j)$,
the Vandermonde determinant.
Thus,
$$
\ZZ[U_n]^{\A_n} = \ZZ[s_1,\dots,s_{n-1}, s_n^{\pm 1}] \oplus 
d \ZZ[s_1,\dots,s_{n-1}, s_n^{\pm 1}] 
$$
This is {\cm}, being free over  $\ZZ[s_1,\dots,s_{n-1}, s_n^{\pm 1}]$.

\subsubsection{$\Sy_n$-lattices} \label{SSS:syn}
If $L$ is a lattice for the symmetric group $\Sy_n$ such that $\ZL^{\Sy_n}$ is {\cm} then 
the Theorem implies that $\Sy_n$ acts as a bireflection group on $L$, and hence on all
simple constituents of  $L \otimes_{\ZZ} {\QQ}$. The simple $\QQ[\Sy_n]$-modules are the
Specht modules $S^\lambda$ for partitions $\lambda$ of $n$. If $n \ge 7$ then the only
partitions $\lambda$ so that $\Sy_n$ acts as a bireflection group on $S^\lambda$ are
$(n)$, $(1^n)$ and $(n-1,1)$; this follows from the lists in \cite{huff} and \cite{wales}.
The corresponding Specht modules are trivial module, $\QQ$, the sign module
$\QQ^-$, and the rational root module $A_{n-1}\otimes_{\ZZ}\QQ$, where $A_{n-1} = 
\{ \sum_i z_ie_i \in U_n \mid \sum_i z_i = 0\}$ and $U_n$ is as in \S\ref{SSS:subgroups}. 
Thus, if $n \ge 7$ and $\ZL^{\Sy_n}$ is {\cm} then we must have
$$
L \otimes_{\ZZ} {\QQ} \cong \QQ^r \oplus \left(\QQ^-\right)^s \oplus 
\left(A_{n-1}\otimes_{\ZZ}\QQ\right)^t
$$
with $s+t \le 2$. In most cases, $\ZL^{\Sy_n}$ is easily seen to be {\cm}. Indeed, we may
assume $r=0$ by Corollary~\ref{SS:rational}. If $s+t \le 1$ then $\Sy_n$ acts as a reflection group
on $L$ and so $\ZL^{\Sy_n}$ is {\cm} by \cite{semigroup}. For $t = 0$ we may quote
the last remark in \S\ref{SS:coefficients}. This leaves the cases $s=t=1$ and $s=0$, $t=2$ to
consider.

If $s=t=1$ then add a copy of $\QQ$ so that $L$ is rationally isomorphic to $U_n \oplus \ZZ^-$.
Using the notation of \S\ref{SSS:subgroups} and putting $t = \x^{(0_{U_n},1)} \in \ZZ[U_n \oplus \ZZ^-]$
the invariants are:
$$
\ZZ[U_n \oplus \ZZ^-]^{\Sy_n} = R \oplus R\varphi
$$
with $R = \ZZ[s_1,\dots,s_{n-1},s_n^{\pm 1}, t + t^{-1}]$ and $\varphi = \frac{1}{2}(\Delta_+ + \Delta)t
+ \frac{1}{2}(\Delta_+ - \Delta)t^{-1}$.

If $s=0$ and $t=2$ then we may replace $L$ by the lattice $U_n^2 = U_n \oplus U_n$. By Lemma~\ref{SS:coefficients} 
$\ZZ[U_n^2]^{\Sy_n}$ is {\cm} precisely if $\FF_p[U_n^2]^{\Sy_n}$ is {\cm} for all primes $p \le n$.
As in \S\ref{SSS:subgroups}, one sees that $\FF_p[U_n^2]^{\Sy_n}$ is a localization of the algebra
``vector invariants" $\FF_p[x_1,\dots,x_n,y_1,\dots,y_n]^{\Sy_n}$. By \cite[Corollary 3.5]{gK01}, 
this algebra is
known to be {\cm} for $n/2 < p \le n$, but the primes $p \le n/2$ apparently remain to be dealt with.

\subsubsection{Ranks $\le 4$} 
As was pointed out in \S\ref{SS:coefficients}, $\ZL^\G$ is always {\cm} when $\rank L \le 2$.

For $L = \ZZ^3$, there are $32$ $\QQ$-classes of finite subgroups $\G \le \GL_3(\ZZ)$. All
$\G$ are solvable; in fact, their orders divide $48$. The Sylow $3$-subgroup $\H \le \G$, if
nontrivial, is generated by a bireflection of order $3$. Thus, $\FF_3[L]^\H$ is
{\cm}, and hence so is $\FF_3[L]^\G$. Therefore, by Lemma~\ref{SS:coefficients}, 
$\ZL^\G$ is {\cm} if and only if $\FF_2[L]^\G$ is {\cm}, and for this to occur, $\G$ must be
generated by bireflections. It turns out that $3$ of the $32$ $\QQ$-classes consist of
non-bireflection groups; these classes are represented by the cyclic groups
$$
\left<\left(\begin{smallmatrix}
-1 & & \\  & -1 & \\ && -1\end{smallmatrix}\right)\right>\ ,
\quad
\left<\left(\begin{smallmatrix}
& 1 & \\  -1 & & \\ && -1\end{smallmatrix}\right)\right>\ ,
\quad
\left<\left(\begin{smallmatrix}
& & -1 \\  -1 & & \\ & -1 & \end{smallmatrix}\right)\right>
$$
of orders $2$, $4$ and $6$ (the latter two classes each split into $2$ $\ZZ$-classes). 
For the $\QQ$-classes consisting of bireflection groups,  Pathak \cite{pkthesis} has checked 
explicitly that  $\FF_2[L]^\G$ is indeed {\cm}.

In rank $4$, there are $227$ $\QQ$-classes of finite subgroups $\G \le \GL_4(\ZZ)$. All but
$5$ of them consist of solvable groups and $4$ of the non-solvable classes are
bireflection groups, the one exception being represented by $\Sy_5$ acting on the signed
root lattice $\ZZ^- \otimes_{\ZZ} A_{4}$. Thus, if the group $\G/\tworef{\G}$ is perfect then
it is actually trivial, that is, $\G$ is a bireflection group. It also turns out that, in this case,
all isotropy groups $\G_m$ are bireflection groups. There are exactly $71$ $\QQ$-classes
that do not consist of bireflection groups. By the foregoing, they lead to non-{\cm}
multiplicative invariant algebras. The $\QQ$-classes consisting of bireflection groups
have not been systematically investigated yet.
The searches in rank $4$ were done with \cite{GAP4}.


\begin{ack}
Some of the research for this article was carried out during a
workshop in Seattle (August 2003) funded by 
Leverhulme Research Interchange Grant F/00158/X
and during the symposium ``Ring Theory" in Warwick, UK (September 2003).
The results described here were reported in the special session
``Algebras and Their Representations" at the AMS-meeting
in Chapel Hill (October 2003). Many thanks to Bob Guralnick for
his helpful comments on an earlier version of this article.
\end{ack}


\end{document}